\documentclass{gtart_h}  

\def\ifplaintex{\expandafter\ifx\csname documentclass\endcsname\relax}

\ifplaintex 
\else
\fi

\expandafter\ifx\csname beginpicture\endcsname\relax
\expandafter\ifx\csname documentclass\endcsname\relax
\input pictex \else
\input prepictex \input pictex \input postpictex \fi\fi

\def\gt{{\mathsurround=0pt\it $\cal G\mskip-2mu$eometry \&\ 
$\cal T\!\!$opology}}        

\def\gtp{{\mathsurround=0pt\it $\cal G\mskip-2mu$eometry \&\ 
$\cal T\!\!$opology $\cal P\!$ublications}}  

\def\lognumber#1{\def\thelognumber{#1}}
\def\volumenumber#1{\def\thevolumenumber{#1}}
\def\papernumber#1{\def\thepapernumber{#1}}
\def\volumeyear#1{\def\thevolumeyear{#1}}

\def\pagenumbers#1#2{\def\startpage{#1}\def\finishpage{#2}}
\def\published#1{\def\publishdate{#1}}
\def\proposed#1{\def\theproposer{#1}}
\def\seconded#1{\def\theseconders{#1}}
\def\received#1{\def\receiveddate{#1}}

\def\accepted#1{\def\accepteddate{#1}}

\long\def\asciiabstract#1{\long\def\theasciiabstract{#1}}

\let\\\par\let\thelognumber\relax
\let\thevolumenumber\relax\let\thepapernumber\relax
\let\thevolumeyear\relax\let\thesamplenumber\relax\let\startpage\relax
\let\finishpage\relax\let\publishdate\relax\let\receiveddate\relax
\let\reviseddate\relax\let\accepteddate\relax\let\theasciititle\relax
\let\theasciiauthors\relax
\let\theasciiabstract\relax
\let\theasciiemail\relax\let\theshortauthors\relax\let\theshorttitle\relax

\long\def\maketitlep{   

\count0=\startpage

\gt\hfill      
\beginpicture
\setcoordinatesystem units <0.33truein, 0.33truein> point at 2.2 0.9
\setplotsymbol ({$\cal G$})
\plotsymbolspacing=9truept
\circulararc 315 degrees from 0 1 center at 0 0
\setplotsymbol ({$\cal T$})
\circulararc 315 degrees from 1 -1 center at 1 0
\endpicture
\break
{\small\ifx\thesamplenumber\relax 
Volume \else Sample
\fi\thevolumenumber\ (\thevolumeyear)
\startpage--\finishpage\nl
Published: \publishdate}
\vglue 0.5truein plus 0.4fil minus 0.1truein

{\parskip=0pt\leftskip 0pt plus 1fil\def\\{\par\smallskip}{\ifplaintex\large
\else\Large\fi\bf\thetitle}\par\medskip}   

\vglue 0pt plus 0.1fil 

{\parskip=0pt\leftskip 0pt plus 1fil\def\\{\par}{\sc\theauthors}
\par\medskip}

\vglue 0pt plus 0.1fil 

{\small\parskip=0pt\let\newline\\
{\leftskip 0pt plus 1fil\def\\{\par}{\sl\theaddress}\par}
\expandafter\ifx\theemail\relax    
\relax\else\vglue 5pt plus 0.02fil minus 2pt\def\\{\stdspace{\rm 
and}\stdspace} 
\cl{Email:\stdspace\tt\theemail}\fi
\ifx\theurl\relax                  
\relax\else\vglue 5pt plus 0.02fil minus 2pt\def\\{\stdspace{\rm 
and}\stdspace}
\cl{URL:\stdspace\tt\theurl}\fi\par}

\vglue 7pt plus 0.3fil minus 3pt

{\bf Abstract}
\vglue 5pt plus 0.1fil minus 2pt

\theabstract

\vglue 7pt plus 0.3fil minus 3pt

{\bf AMS Classification numbers}\quad Primary:\quad \theprimaryclass

Secondary:\quad \thesecondaryclass

\vglue 5pt plus 0.3fil minus 2pt

{\bf Keywords:}\quad \thekeywords

\vglue 10pt plus 0.5fil minus 5pt

{\small  Proposed: \theproposer\hfill Received: \receiveddate\nl
Seconded: \theseconders\hfill 
\ifx\reviseddate\relax                         
Accepted: \accepteddate                        
\else
Revised: \reviseddate                          
\fi}
\eject
}       

\let\maketitlepage\maketitlep
\let\maketitle\maketitlepage

\font\phead=cmsl9 scaled 950
\font=cmsl9 scaled 1050
\font\pnum=cmbx10 scaled 913
\font\lnum=cmbx10 
\font\pfoot=cmsl9 scaled 950
\font=cmsl9 scaled 1050
\ifplaintex
\headline{\vbox to 0pt{\vskip -4.5mm\line{\small\phead\ifnum
\count0=\startpage ISSN 1364-0380 (on line)
1465-3060 (printed) \hfill {\pnum\folio}\else\ifodd\count0\def\\{ }%
\ifx\theshorttitle\relax\thetitle\else\theshorttitle\fi\hfill{\pnum\folio}
\else\def\\{ and }{\pnum\folio}\hfill\ifx\theshortauthors\relax\theauthors
\else\theshortauthors\fi\fi\fi}\vss}}
\footline{\vbox to 0pt{\vglue 0mm\line{\small\pfoot\ifnum\count0=\startpage
\copyright\ \gtp\hfill\else
\gt, Volume \thevolumenumber\ (\thevolumeyear)\hfill\fi}\vss
}}
\else
\makeatletter
\def\@oddhead{{\small\lhead\ifnum\count0=\startpage ISSN 1364-0380 (on line)
1465-3060 (printed) \hfill {\lnum\number\count0}\else\ifodd\count0
\def\\{ }\ifx\theshorttitle\relax \thetitle \else\theshorttitle\fi\hfill
{\lnum\number\count0}\else\def\\{ and }{\lnum\number\count0}
\hfill\ifx\theshortauthors\relax 
\theauthors\else\theshortauthors\fi\fi\fi}}\def\@evenhead{\@oddhead}
\def\@oddfoot{\small\lfoot\ifnum\count0=\startpage\copyright\ \gtp\hfill\else
\gt, Volume \thevolumenumber\ (\thevolumeyear)\hfill\fi}
\def\@evenfoot{\@oddfoot}
\makeatother
\fi

\newwrite\gtoutfile
\long\gdef\makeheadfile{  
{\def\\{, }\def\s{ }
\immediate\openout\gtoutfile head.xxx
\immediate\write\gtoutfile{Proxy-for: \ifx\theasciiauthors\relax
\theauthors\else\theasciiauthors\fi\s<\ifx\theasciiemail\relax\theemail\else\theasciiemail\fi>}
\immediate\write\gtoutfile{\noexpand\\}
\immediate\write\gtoutfile{Authors: \ifx\theasciiauthors\relax
\theauthors\else\theasciiauthors\fi}
{\def\\{ }\immediate\write\gtoutfile{Title: \ifx\theasciititle\relax
\thetitle\else\theasciititle\fi}}
\immediate\write\gtoutfile{Subj-class: GT or SG or MG etc}
\immediate\write\gtoutfile{MSC-class: \theprimaryclass\ifx\thesecondaryclass\relax\else, \thesecondaryclass\fi}
\immediate\write\gtoutfile{Journal-ref: Geom. Topol. \thevolumenumber
(\thevolumeyear) \startpage-\finishpage}
\immediate\write\gtoutfile{Comments: Published by Geometry and Topology at}
\immediate\write\gtoutfile{\s\s http://www.maths.warwick.ac.uk/gt/GTVol\thevolumenumber/paper\thepapernumber.abs.html}
\immediate\write\gtoutfile{\noexpand\\}
\immediate\write\gtoutfile{}
\ifx\theasciiabstract\relax
\immediate\write\gtoutfile{\theabstract}\else
\immediate\write\gtoutfile{\theasciiabstract}\fi
\immediate\write\gtoutfile{}
\immediate\write\gtoutfile{\noexpand\\}
\immediate\write\gtoutfile{}
\immediate\closeout\gtoutfile}}  

\def\maketitlepage{\maketitlep\makeheadfile}
\let\maketitle\maketitlepage

\lognumber{454}
\volumenumber{8}\papernumber{27}\volumeyear{2004}
\pagenumbers{1013}{1031}
\received{13 May 2004}
\published{7 August 2004}
\accepted{11 July 2004}
\proposed{Peter Ozsvath}
\seconded{Tomasz Mrowka, Peter Kronheimer}

\theoremstyle{plain}
\usepackage{amsmath, amssymb,amscd}

\newtheorem{lem}{Lemma}[section]
\newtheorem{prop}[lem]{Proposition}
\newtheorem{thrm}[lem]{Theorem}
\newtheorem{cor}[lem]{Corollary}

\newtheorem{tm}{Theorem}

\newcommand{\Ozsvath}{{Ozsv{\'a}th  }}
\newcommand{\Szabo}{{Szab{\'o} }}

\newcommand{\Z}{{\bf{Z}}}
\newcommand{\Q}{{\bf{Q}}}
\newcommand{\coker}{{{\text {coker} \ }}}
\newcommand{\rank}{{{\text {rank} \ }}}
\newcommand{\ts}{{{\thinspace}}}
\newcommand{\gr}{{{\mathrm{gr} }}}

\newcommand{\hfp}{{{HF^+}}}

\newcommand{\hfred}{{{HF^{\rm red}}}}

\newcommand{\hfi}{{{HF^{\infty} }}}

\newcommand{\spi}{{{\mathfrak s}}}
\newcommand{\vspi}{{{\mathfrak t}}}
\newcommand{\uspi}{{{\mathfrak u}}}
\newcommand{\vvspi}{{{\mathfrak v}}}
\newcommand{\spinc}{{{\mathrm{Spin}^c}}}

\begin{document}
\title{Lens space surgeries and a conjecture of\\Goda and Teragaito}
\author{Jacob Rasmussen}
\address{Department of Mathematics, Princeton University\\Princeton, 
NJ 08544, USA}
\email{jrasmus@math.princeton.edu}

\begin{abstract}
Using work of \Ozsvath and Szab{\'o}, we show that if a nontrivial
knot in \(S^3\) admits
a lens space surgery with slope \(p\), then \(p \leq 4g+3\), where \(g\)
is the genus of the knot. This is a close approximation to a bound
conjectured by Goda and Teragaito. 
\end{abstract}

\asciiabstract{Using work of Ozsvath and Szabo, we show that if a
nontrivial knot in S^3 admits a lens space surgery with slope p, then
p <= 4g+3, where g is the genus of the knot. This is a close
approximation to a bound conjectured by Goda and Teragaito.}

\primaryclass{57M25}
\secondaryclass{57R58}
\keywords{Lens space surgery, Seifert genus, Heegaard Floer homology}

\maketitlepage

\section{Introduction}

Let \(K\) be a knot in \(S^3\), and denote by \(K_r\) the
three-manifold obtained by performing Dehn surgery on \(K\) with slope
\(r=p/q\). If \(K_r\) is homeomorphic to a lens space we say that \(K\)
 admits a lens space surgery with slope \(r\).
 In recent years, Kronheimer, Mrowka, Ozsv{\'a}th, and \Szabo have
used Floer homology for three-manifolds to give constraints on such
knots \cite{OS4}, \cite{OSLens}, \cite{KMOS}. Generally speaking, these
constraints are derived from the fact that lens spaces belong to a larger
class of spaces, known as {\it L--spaces}, for which the reduced Floer
homology groups \(\hfred\) vanish. 

On the other hand there are many \(L\)--spaces which are not lens
spaces. In particular, if \(K\) admits a single \(L\)--space surgery with
positive slope, then \(K_p\) is an \(L\)--space for every integer  
 \(p \geq 2g(K)-1 \), where \(g(K) \) denotes the genus of \(K\)
 \cite{KMOS}. In contrast, when \(K\) is hyperbolic,
 the cyclic surgery theorem of \cite{CyclicS} tells us that at
 most two of these surgeries are actually lens spaces. In this
note, we show that Floer homology can be used to distinguish at least some of
these \(L\)--space surgeries from lens spaces. In particular, we
prove the following result.
\begin{tm}
Suppose \(K\) is a nontrivial knot which
 admits a lens space surgery of slope \(r\). Then
\begin{equation*}
 |r| \leq 4 g(K) + 3. 
\end{equation*}
\end{tm}
\noindent The inequality is sharp --- equality holds for the case of \(4k+3\)
surgery on the right-handed \((2,2k+1)\) torus knot, which gives the lens space
\(L(4k+3,4)^*\). 

This result closely approximates a bound conjectured by Goda
and Teragaito in \cite{GT}. More specifically, they showed that if
\(K\) is a {\it hyperbolic} knot which admits a lens space surgery of
slope \(p\), then \(|p| \leq 12g(K)-7 \), and conjectured that in fact
\begin{equation*}
2g(K)+8\leq |p| \leq 4g(K)-1. 
\end{equation*}
Something close to the first inequality was proved in
Corollary 8.5 of
\cite{KMOS}, where it was shown that if \(K\) admits an \(L\)--space
surgery of slope \(p\), then \( 2g(K)-1 \leq |p| \). 
Theorem 1 seems to be a natural (and only minimally weaker)
reformulation of the second inequality which applies to all
knots.

The proof of the theorem is based on work of \Ozsvath and \Szabo in
\cite{OS4}. In addition, we use an analog of an inequality of
Fr{\o}yshov \cite{Froy3} and some elementary facts about Dedekind
sums. The paper is arranged as follows: in section 2, we review the
results of \cite{OS4} and outline the proof of Theorem 1. Section 3
is devoted to the proof of Fr{\o}yshov's inequality, and section 4
contains the necessary results on Dedekind sums.

Throughout this note, we work in the category of oriented
manifolds. All maps, homeomorphisms, {\it etc.} are assumed to be
orientation preserving unless specified otherwise. For lens spaces, 
our orientation convention is the one used in 
\cite{GompfStip} and \cite{Rolfsen}, namely, that 
\(-p\) surgery on the unknot produces the oriented lens space
\(L(p,1)\). (Note that this is the opposite of the convention used in
\cite{OS4} and \cite{KMOS}).

\medskip
\noindent{\bf Acknowledgements}\qua The author would like to thank Peter
\Ozsvath and Zoltan \Szabo for helpful conversations. This work was
partially supported by an NSF Postdoctoral Fellowship. 

\section{Outline of proof}
\label{Sec:Proof}
Suppose that \(K\) is a nontrivial knot and that \(K_r\) is a lens
 space. Without loss of generality, we may assume that \(r\) is an
 integer. Indeed,
 the cyclic surgery theorem implies that this must be case unless
 \(K\) is a torus knot. On the
 other hand, if \(K\) is the right-handed \((a,b)\) torus knot, it is well
 known \cite{Moser} that \(K_{p/q}\) is a lens space if and only if
\(qab - p = \pm 1\). In particular, the slope  attains its
 largest value when \(p/q = ab + 1\) is an integer.

We now review the results of \cite{OS4} on knots  admitting integral
 lens space surgeries. For technical reasons, it is convenient to
 assume that the slope of these surgeries is negative. By considering the
mirror image of \(K\), if necessary, we may arrange that this is the
 case. From this point on, then, we will assume that \(K_{-p} \) is a lens
 space, where \(p\) is a positive integer. 

\subsection{The exact triangle} 
\label{SubSec:Triangle} 
 Let \(W_1\) be the surgery cobordism from \(K_0\) to
\(S^3\), and let \(x\) be a generator of \(H^2(W_1) \cong \Z\). We use
 the notation \( \spi_i\) to refer either to the \(\spinc\) structure
 on \(W_1\) with \(c_1(\spi_i) = 2ix\) or its restriction to
 \(K_0\). (It should be clear from context which manifold is being
 considered.) Likewise, if  \(W_2\) is the surgery cobordism from \(S^3 \) to
\(K_{-p}\) and \(y \in H^2(W_2)\) is a generator, we let \(\vspi_i\)
 be the \(\spinc\) structure on \(W_2\) with \(c_1(\vspi_i) =
 (-p+2i)y\), and \(\vspi_i '\) be its restriction to \(K_{-p}\). Note
 that 
\(\vspi_i '\) only depends on the value of \(i \ \text{mod} \
 p\). 

 The exact triangle with twisted coefficients
 \cite{OS2} gives a long exact sequence
\begin{equation*}
\begin{CD}
 \oplus_{i\equiv k \ts (p)} \underline{HF}^+ (K_0, \spi_i) \!
@>{\oplus F_{W_1,\spi_i}}>>\!
 \hfp (S^3)[T,T^{-1}]\! @>{ F_{W_2,k}}>>
\hfp (K_{-p}, \vspi_k')[T,T^{-1}] 
\end{CD}
\end{equation*}
where
\begin{equation*}
 F_{W_2,k} (x) = \sum _{n \in \Z} T^{n}\cdot F_{W_2,\vspi_{k+pn}}(x).
\end{equation*}
Let \(h_i\) be the rank of 
 \(F_{W_1,\spi_i} \co \underline{HF}^+
 (K_0,\spi_i) \to \hfp(S^3)[T,T^{-1}] \), viewed as a map of \(\Z[T,T^{-1}]\)
 modules. Combining 
 results from \cite{OS4} and \cite{OSGen} gives the following:

\begin{prop}
Suppose \(K_{-p}\) is a lens space. Then \(p \geq 2g(K) - 1\), and
\(h_i\) is nonzero if and only if \( -g(K)+1 \leq  i \leq g(K) -1\). 
Moreover, for \(-p/2 \leq k \leq p/2 \) \em{
\begin{equation*}
 \rank \ker F_{W_2,k} = \rank  \underline{HF}^+ (K_0, \spi_k) = h_k.
\end{equation*}}
\end{prop}

\begin{proof}
Since \(K_{-p}\) is a lens space, \(\hfp(K_{-p}, \vspi_k') \cong \hfp (S^3) 
\cong \Z[u^{-1}]\). Now \(\underline{HF}^+(K_0, \spi_i)\) has finite
rank as a \( \Z[T, T^{-1}]\) module, and  \(F_{W_2,k} \) is  a
\(u\)--equivariant map. From this, it follows that  
 \( \ker F_{W_2,k}
\) must be of the form \(\langle 1, u^{-1},\ldots, u^{-n}\rangle \otimes
\Z[T,T^{-1}]\) for some value of \(n\), and that \( \coker F_{W_2,k}
\) must have zero rank. Thus  the maps \(
F_{W_1,\spi_i} \) must all be of full rank, and the sum 
\begin{equation}
\label{Eq:Sum}
 \bigoplus_{i\equiv k \ts (p)} \underline{HF}^+ (K_0, \spi_i)
\end{equation}
can contain at most one term of nonzero rank.
In this case, we have
\begin{equation*}
 \rank \ker F_{W_2,k} = \rank  \underline{HF}^+ (K_0, \spi_i) = h_i.
\end{equation*}
where \(i\) is the index of the nontrivial summand. 
Since \(h_i \leq
h_j\) whenever \(|i|>|j|\) (Proposition 7.6 of \cite{thesis}), it
follows that \(i\) must be the representative of \(k \ \text{mod} \ p\)
with smallest absolute value, so \(-p/2 \leq i \leq p/2\).

According to \cite{OSGen}, the largest value of \(i\) for which
\(\hfp(K_0, \spi_i)\) is nontrivial is \(g(K)-1\). 
It follows that \(h_i\) is nonzero if and only if 
\begin{equation*}
 -g(K)+1 \leq i \leq g(K) -1.
\end{equation*} 
Thus there are \(2g(K)-1\) values of \(i\) for which \(h_i\)
is nontrivial, and the condition that the sum in equation
(\ref{Eq:Sum}) contains at most one
nontrivial term for all values of \(k\) is equivalent to the statement that
\(p \geq 2g(K) - 1\). 
\end{proof}

\subsection{The \(d\)--invariant}
Let \(Y\) be a rational homology three-sphere, and let \( \spi \) be a
 \( \spinc \) structure on it.
Then \(\hfp(Y, \spi) \otimes \Q \) is absolutely graded and 
contains a unique summand isomorphic to \(\Q[u^{-1}] \). The invariant \(d(Y,
\spi) \) is defined \cite{OS4}
 to be the absolute grading of \( 1 \in \Q[u^{-1}]
\). (This is the analog of Fr{\o}yshov's \(h\) invariant in
Seiberg--Witten theory.)  

The \(d\)--invariants of \(K_{-p}\) are easily expressed in terms of
the \(h_k\)'s. Recall that the map \(F_{W_2,k}\) is a sum of maps
\( F_{W_2,\vspi_i}  \co \hfp (S^3) \to \hfp (K_{-p}, \vspi_k') \). Since the
intersection form on \(W_2\) is negative-definite, each \(
F_{W_2,\vspi_i}\) is a surjection. Moreover, \(F_{W_2,\vspi_i} \) is 
 a graded map; it shifts the absolute
grading by a fixed rational number \(\epsilon(p,i)\) which depends
 only on homological data associated to the cobordism \(W_2\). In
 particular, this number is independent of \(K\).

The kernel of \(F_{W_2,k}\) is generated by \(\langle 1,
u^{-1}, \ldots, u^{-h_k+1}\rangle\), so  \(u^{-h_k} \in
\hfp(S^3) \) must map to a nontrivial multiple of 
\(1 \in \hfp (K_{-p}, \vspi_k)[T,T^{-1}]\). 
Thus if we let 
\begin{equation*}
 E(p,k) = \max \ts \{ \epsilon(p,i) \ts | \ts i \equiv k \ts
(p) \}
\end{equation*}
 it follows that for \(-p/2 \leq k \leq p/2 \), 
\begin{equation}
\label{Eq:DRelation}
d (K_{-p}, \vspi_k') = E(p,k) + 2h_k. 
\end{equation}
Specializing to the case where \(K\) is the unknot, we
see that the set of \(d\)--invariants of the lens space \(L(p,1) \) is 
\(\{E(p,k) \ts | \ts k \in \Z/p\} \).

\subsection{The Casson--Walker invariant}
If \(Y\) is an integral homology three-sphere, then
 its Casson invariant \( \lambda (Y) \) is determined
by \(d(Y)\) and the group \( \hfred (Y)\) (Theorem 5.1 of
 \cite{OS4}.) More generally, if \(Y\) is a rational homology
 three-sphere, it is expected that its
 Casson--Walker invariant \( \lambda (Y)\) 
will be related to  \( \hfred (Y)\) and the \(d\)--invariants of \(Y\). 
This relation is particularly simple when \(Y\) is a lens space:
\begin{lem}
\label{Lem:DCW}
\begin{equation*}
  \sum_{\spi} d(L(p,q), \spi) =  p \lambda (L(p,q)))
\end{equation*}
where the sum runs over all \( \spinc \) structures on \(L(p,q)\).
\end{lem}
\noindent The proof will be given in section~\ref{Sec:Dedekind}.
Combining  with equation (\ref{Eq:DRelation}), we get
\begin{equation*}
p(\lambda (K_{-p}) - \lambda (L(p,1))) = 2 \sum _{i=-g+1}^{g-1} h_i.
\end{equation*}

\subsection{Fr{\o}yshov's inequality}
The new geometric input in the proof of Theorem 1 is the
following fact, which is analogous to a theorem of
Fr{\o}yshov in instanton Floer homology \cite{Froy3}. Its proof is the
subject of section~\ref{Sec:Froyshov}.

\begin{thrm}
\label{Prop:Froyshov}
 Let \(K\) be a knot in \(S^3\), and let
  \(g_*(K) \) be its slice genus. Then \(h_i(K) = 0\) for \(|i|\geq
  g_*(K)\), while for \(|i|< g_*(K) \)
\begin{equation*}
h_i(K) \leq \Bigr\lceil \frac{g_*(K) - |i|}{2} \Bigl\rceil.
\end{equation*}
\end{thrm}
\noindent If \(K\) admits a lens space surgery, the results of 
\cite{OSLens} and \cite{OSGen} show that  \(g_*(K)= g(K) \), so
\begin{equation*}
p(\lambda (K_{-p}) - \lambda (L(p,1))) \leq\, 2\!\!\!\! \sum _{i=-g+1}^{g-1} 
\Bigr\lceil \frac{g_*(K) - |i|}{2} \Bigl\rceil = 
 g(K)(g(K)+1).
\end{equation*}

\subsection{Proof of the theorem}
Suppose that the inequality is false. Then \( p \geq 4g(K) +4\), so 
\(g(K)+1 \leq \frac{p}{4} \). Substituting into the previous
inequality,
 we find that 
\begin{equation*}
\lambda (K_{-p}) - \lambda (L(p,1)) \leq \frac{1}{4}(\frac{p}{4}-1).
\end{equation*}
The value of \( \lambda (L(p,q)) \) is given by a certain arithmetic function of \(p\) and \(q\), known as a {Dedekind sum}. The following purely arithmetic result is proved in  section~\ref{Sec:Dedekind}:
\begin{prop}
\label{Prop:DBound}
Suppose that \(Y\) is a lens space with \( |H_1(Y)| = p\), and that 
\begin{equation*}
\lambda (Y) - \lambda (L(p,1)) \leq \frac{1}{4}(\frac{p}{4}-1).
\end{equation*}
Then \(Y\) is homeomorphic to one of \(L(p,1)\), \(L(p,2)\), or  \(L(p,3)\). 
\end{prop}

The possibility that \(K_{-p}\) is homeomorphic to
 \(L(p,1)\) is ruled out by the main theorem of \cite{KMOS}. To
 eliminate the other two cases, we use the following proposition, which is
 also proved in section~\ref{Sec:Dedekind}.

\begin{prop}
\label{Prop:23}
If \(K_{-p}\) is homeomorphic to \(L(p,2)\), then \(p=7\). If
it is homeomorphic to \( L(p,3)\), then either \(p=11\) or \(p=13\). 
\end{prop}

From the
table at the end of \cite{OS4}, we see that if \(L(7,2)\), \(L(11,3)\),
or \(L(13,3)\) is 
realized by integer surgery on a knot \(K\), then \(K\) must have genus
1, 2,  or 3, respectively, and the inequality of
Theorem 1 is satisfied.
(In fact, these lens spaces are realized by surgery on the
torus knots \(T(2,3)\), \(T(2,5)\), and \(T(3,4)\), respectively.) 
 This concludes the proof of the theorem. \qed

\section{Fr{\o}yshov's inequality}
\label{Sec:Froyshov}
We now turn to the proof of Theorem~\ref{Prop:Froyshov}. 
 The argument we give 
 is essentially that of \cite{Froy3}, but adapted along the
lines of \cite{OS4} to fit the
Heegaard Floer homology. We begin by reformulating the problem
slightly. Let \(K\) be a knot in \(S^3\), and choose \(n \gg 0\). 
Then for \(-n/2 \leq k \leq n/2\), we have
 \(d(K_{-n},\vspi_k') = E(n,k) + 2h_k (K).\) 
To prove the theorem, we will estimate the size of 
\( d(K_{-n},\vspi_k') \). 

To this end, we consider
the surgery cobordism  \({W_2}^*\) from \(K_{-n}\) to \(S^3 \).
(This is  the cobordism \(W_2\) of section~\ref{SubSec:Triangle} 
with its orientation reversed.)
We fill in the \(S^3\)
 boundary component of \({W_2}^*\) with a four-ball to get a
  four-manifold \(W'\). Then 
 \(H_2(W') \cong \Z \), and the generator of this group can be
 represented by an embedded surface \( \Sigma_g \) with genus \(g = g_*(K)
 \) and self-intersection \(n\). Finally, let \(W\) be the
 four-manifold obtained by removing a tubular neighborhood  of \(\Sigma_g
 \) from \(W'\). \(W\) is a cobordism from \(K_{-n}\) to
 the circle bundle over \(\Sigma_g\) with Euler number \(-n\), which we
 denote by \(B_{-n}\). This choice of name is a natural one, since
 \(B_{-n}\) can be obtained by doing
 \(-n\) surgery on the ``Borromean
 knot'' \(B \subset \#^{2g} (S^1\times S^2) \). 

We now consider the topology of the cobordism \(W\).
An easy computation shows that \(H^2(W) \cong H^2(B_{-n}) \cong 
\Z^{2g} \oplus \Z/n \), and the restriction map to
\(H^2({K^*_{-n}})\) is projection onto the second factor. It
follows that there is a unique torsion \( \spinc \) structure on \(W\)
which restricts to \( \vspi_k'\) on \({K^*_{-n}} \). We
denote this \(\spinc\) structure and its restriction to \(B_{-n}\) by
\( \vspi_k'\) as well.

Note that there is another natural way to label the torsion \( \spinc \)
structures on \(B_{-n}\). Namely, we can view \(B_{-n} \) as \({-n}\)
surgery on the knot \(B\) and use the labeling convention of
section~\ref{SubSec:Triangle}. To be precise, let 
 \(X_2\) be the surgery cobordism from
 \( \#^{2g} (S^1\times S^2) \) to \(B_{-n} \). Then the restriction
 map \(H^2(X_2) \to H^2(\#^{2g} (S^1\times S^2))\) has kernel isomorphic
 to \(\Z\).
 If \(x\) is a generator of this group, we let \(\uspi_k\)
 be the \(\spinc \) structure on \(X_2\) with \(c_1(\uspi_k) = (-n + 2k)x
 \), and \( \uspi_k'\) be its restriction to \(B_{-n} \). 

\begin{lem}
For an appropriate choice of the generator \(x\), we have
\( \vspi_k ' = \uspi_k' \).
\end{lem}

\begin{proof}
Let \(X'\) be the double of \(W'\), and let \(S \subset
X' \) be the embedded sphere which is obtained by gluing together the
cocore of the two-handle in \(W'\) (which is an embedded disk
generating \(H_2(W', \partial W')\)) and its mirror image in
\((W')^*\). \(S\) intersects \( \Sigma_g\) geometrically once. If we
remove a tubular neighborhood of \( \Sigma_g\) from \(X\), we get a
four-manifold \(X\) which is the union of \((W')^*\) and \(W\). \(S
\cap X\) is an embedded disk \(D\) whose boundary 
 is a fiber of the circle bundle \(B_{-n}\). Let
\(D_1\) be the disk \(D \cap (W')^*\). Then 
\(D_1 \) is Poincare dual to the generator  of \(H^2((W')^*) \).

Recall that 
the \(\spinc\) structure \(\vspi_k\) on \(W_2\) was defined by
\(c_1(\vspi_k) = (-n+2k)y \), where \(y\) was a generator of
\(H^2(W_2)\). \((W')^*  = W_2 \cup D^4\), so \(\vspi_k\)
extends uniquely to a \( \spinc \) structure \(\vspi_k\) on \((W')^*\) with
  \( c_1(\vspi_k) = (-n + 2k) \text{PD}(D_1) \). Now let \(\vvspi_k \)
be the \( \spinc \) structure on \(X\) with \( c_1(\vvspi_k) =
(-n+2k)\text{PD}(D) \). Then \(\vvspi_k|_{(W')^*} = \vspi_k\), and
\(\vvspi_k|_{B_{-n}}\) is torsion,  so
\(\vvspi_k|_{W} = \vspi_k'\).  On the other hand, if
\(X'' \subset X\) is a regular neighborhood of
 \(B_{-n} \cup D\), it is not difficult to see that
 \(X'' \cong X_2\), and that the kernel of the restriction map
\(H^2 (X'') \to H^2(\#^{2g} (S^1\times
 S^2))\)  is generated by \( \text{PD}(D)\). Thus
\(\vvspi_k |_{X''} = \uspi_k \), and the claim follows. 
\end{proof}

Returning to the topology of \(W\), we further calculate
 that \(b_1(W) = b_2^+(W) =
0\). Now if  \(W\) were a cobordism between two rational homology
spheres, the fact that  \(b_2^+(W) = 0 \) would imply that the induced
map \(F_{W,\spi}^\infty\) is an isomorphism for any  \(
\spinc \) structure \( \spi \) on \(W\). 
In our case, \(B_{-n}\) is not a homology
sphere, so the situation is somewhat more complicated. Nevertheless,
it is still true that \(F_{W,\vspi_k'}^\infty\) is an injection:

\begin{lem}
Suppose \(W\) is a cobordism from \(Y_1\) to \(Y_2\) and that
\begin{equation*}
b_1(Y_1) = b_1(W) = b_2^+(W) = 0.
\end{equation*}
Let \(\spi\) be a \( \spinc \) structure on \(W\) whose restriction
\(\spi_i\) to \(Y_i\) is torsion, and suppose moreover that
 \(\hfi(Y_2, \spi_2)\) is ``standard,'' in the sense
that its rank as a \(\Z[u,u^{-1}]\) module is \(2^{b_1(Y_2)}\). Then 
\begin{equation*} 
F^\infty_{W,\spi} : \hfi (Y_1,\spi_1) \to \hfi
(Y_2, \spi_2) 
\end{equation*}
 maps \( \hfi (Y_1,\spi_1)\) isomorphically onto \(A_{\spi_2}
\subset \hfi (Y_2, \spi_2) \), where \em{
\begin{equation*}
A_{\spi_2} = \{ x \in \hfi(Y_2, {\spi_2}) \ts | \ts \gamma \cdot x = 0 \ \
\forall  \ \gamma \in H_1(Y_2)/\text{Tors} \}.
\end{equation*}}
\end{lem}

\noindent {\bf Remark}\qua If we wish to avoid the use of twisted
coefficients, the condition that \(\hfi (Y_2, \spi_2) \) be standard is
clearly necessary. For example, let \(W\) be the cobordism from
\(S^3\) to \(T^3\) obtained by removing a ball and a neighborhood of a
regular fiber from the rational elliptic surface \(E(1)\). Then \(W\)
satisfies the homological conditions of the lemma, but
\(F^\infty_{W,\spi}\) is the zero map. 

\begin{proof}
This is essentially contained in the proof of Theorem~9.1 in
\cite{OS4}. The argument may be summarized as follows.
 The cobordism
\(W\) can be broken into a composition of three cobordisms
 \(W_i \ (1\leq i \leq 3)\), each corresponding to the addition of
 handles of index \(i\). The hypothesis that \(b_1(W) = b_2^+(W)=0 \)
 implies that \(W_2\) can be further decomposed into a composition of
 cobordisms \(W_{2}^-\), \(W_{2}^0\) and \(W_{2}^+\), where each
 two-handle addition in \(W_{2}^-\) decreases \(b_1\) of the terminal
 end by \(1\), each two-handle addition in \(W_{2}^0\) does not change
 \(b_1\) of the terminal end, and each two-handle in \(W_{2}^+\)
 increases \(b_1\) of the terminal end by~\(1\). 

We further subdivide \(W_{2}^+\) into a sequence of elementary
cobordisms, each corresponding to the addition of a two-handle. Let
\(Y^i\) be the terminal end of the \(i\)-th such cobordism, and let
\(Y^0\) be the initial end of the first one. Then the hypothesis that
\(b_1(W) = 0 \) implies that \(Y^0\) is a rational homology sphere, so
\(b_1 (Y^i) = i\). Likewise, the fact that \(b_2^+(W) = 0 \)
implies that the restriction of \( \spi\) to \(Y^i\) (which we
continue to denote by \(\spi\)) must be torsion. Finally,  the
exact triangle shows that  
\begin{equation*}
\rank \hfi(Y^{i+1}, \spi) \leq 2 \ \rank \hfi
(Y^i, \spi ).
\end{equation*}
 It follows that if \( \hfi (Y^i, \spi) \) is
standard, 
then \(\hfi
(Y^j, \spi)\) is standard for all \(j\leq i\). Since \(W_3\) is composed
entirely of three-handles, the terminal end of \(W_{2}^+\) is
homeomorphic to \(Y_2 \#^n(S^1\times S^2)\). Now \(\hfi (Y_2, \spi_2) \) is
standard by hypothesis, and this implies that \(\hfi (Y_2
\#^n(S^1\times S^2), \spi) \) is standard as well. Thus \(\hfi
(Y^i, \spi)\) is standard for all \(i\). 

One can now  check directly, using Proposition 9.3 of \cite{OS4}, that
\(F^\infty_{W_1\cup W_{2}^-, \spi} \) and \(F^\infty_{W_{2}^0,\spi} \)
 are isomorphisms, that
\(F^\infty_{W_{2}^+, \spi}\) is injective and maps onto \(A_{\spi}\), and
that \(F^\infty_{W_3, \spi} \) preserves this property. 
\end{proof}

In order to apply the lemma, we must check that \( \hfi (B_{-n},
\vspi_k')\) is standard. From the exact triangle
for the knot \(B \subset \#^{2g} (S^1 \times S^2)\), we see that 
\begin{equation*}
\hfi (B_{-n}, \vspi_k') \cong \hfi (\#^{2g} (S^1 \times S^2) , \vspi)
\end{equation*}
where \(\vspi \) is the unique torsion \( \spinc \) structure on 
\( \#^{2g} (S^1 \times S^2) \). The latter group
 is standard, so \( \hfi (B_{-n},\vspi_k')\) must be
 standard as well.

Let \(A_{\vspi_k'} \subset \hfi(B_{-n},\vspi_k') \) be as in the lemma. Then 
\(A_{\vspi_k'} \otimes \Q \cong \Q[u,u^{-1}] \), and its image under the map
\(\pi \co \hfi(B_{-n},\vspi_k') \otimes \Q \to \hfp(B_{-n},\vspi_k')
\otimes \Q \) will be isomorphic
to \(\Q[u^{-1}]\). In analogy with the \(d\)--invariant for rational
homology spheres, we define \(d(B_{-n}, \vspi_k') \) to be the absolute
grading of \( 1 \in \pi(A_{\vspi_k'}) \otimes \Q \cong \Q[u^{-1}] \). 

\begin{lem}
\( d(K_{-n}, \vspi_k') \leq d(B_{-n}, \vspi_k') + g\). 
\end{lem}

\proof
Consider the map \( F^+_{W, \vspi_k'} \co \hfp (K_{-n}, \vspi_k') \to \hfp
(B_{-n}, \vspi_k') \). This map is \(u\)--equivariant and agrees with
\(F^\infty_{W, \vspi_k'} \) in high degrees, which implies that it
takes \(
\pi(\hfi (K_{-n}, \vspi_k')) \) onto  \( \pi(A_{\vspi_k'}) \). Thus if \(1\in 
\pi(\hfi(K_{-n}, \vspi_k')) \otimes \Q \cong \Q[u^{-1}] \)
 is the element with the lowest
absolute grading, we must have
\begin{equation*}
\gr ( F^+_{W, \vspi_k'}(1)) \leq d(B_{-n}, \vspi_k').
\end{equation*}
But \( \gr (1) = d(K_{-n},\vspi_k') \) and \(  F^+_{W, \vspi_k'} \) shifts
the absolute grading by 
$$
\frac{c_1^2(\vspi_k') - 2 \chi (W) - 3 \sigma (W)}{4} = 
\frac{0 - 2\cdot 2g - 3 \cdot 0}{4} = -g.
\eqno{\qed}$$

Since \( d(K_{-n},\vspi_k') = E(n,k) + 2 h_k \), the proof of
Theorem~\ref{Prop:Froyshov}  reduces to the following computation:

\begin{prop}
\label{Prop:dCalc}
For \(n \gg 0\),
\(\displaystyle{
d(B_{-n}, \vspi_k') = E(n,k) - g + 2 \Bigl\lceil \frac{g-|k|}{2} \Bigr\rceil}. \)
\end{prop}

\begin{proof} The Floer homology of \(B_{-n}\) was computed by
  \Ozsvath and \Szabo in section 9 of \cite{OS7}. 
More specifically, they show that the knot Floer homology of 
the Borromean knot  \(B \subset \#^{2g}(S^1\times S^2) \)
 is given by 
\begin{equation*}
\widehat{HFK}(B,i) \cong \Lambda^{g+i}(H^1(\Sigma_g)).
\end{equation*}
This complex is {\it perfect}, in the sense that the
homological grading is equal to the Alexander grading, and there are
no differentials, even in the larger complex \( HFK^\infty (B)
 \cong \widehat{HFK}(B) \otimes \Z[u, u^{-1}] \). \Ozsvath and \Szabo
 also compute the action of
\begin{equation*}
 H_1(\#^{2g} (S^1 \times S^2)) \cong H_1(B_{-n})/\text{Tors} \cong
 H_1(\Sigma_g)  
\end{equation*}
 on \( HFK^\infty (B) \); it is given by
\begin{equation*}
\gamma \cdot (\omega \otimes u^{n}) = \iota_\gamma \omega \otimes u^n
+ \text{PD}(\gamma) \wedge \omega \otimes u^{n+1},
\end{equation*}
where \(\text{PD}(\gamma)\) denotes the Poincare dual of \(\gamma \)
viewed as an element of \(H_1(\Sigma_g) \). 

Let \(\{a_1,\ldots, a_g, b_1, \ldots, b_g \}\) be a symplectic basis of 
\(H^1(\Sigma_g)\). We can write any
 \( \omega \in {HFK}^\infty (B)\) in the form
 \( \omega = \omega_1 + b_1 \wedge \omega_2 \), where \(b_1\)
does not appear in the expressions for \( \omega_1\) and \( \omega_2
\). Then 
\begin{equation*}
 \text{PD}(a_1) \cdot (\omega _1 + b_1 \wedge \omega _2) 
= \omega_2 + a_1 \wedge \omega_1 \otimes u +a_1 \wedge b_1 \wedge
\omega_2 \otimes u.
\end{equation*}
For this expression to vanish, we must have \( \omega_2 = -a_1 \wedge
 \omega _1 \otimes u \). Thus \( \omega = (1 + a_1 \wedge b_1 \otimes
 u) \omega_1\), where \(b_1\) does not appear  in the expression for  \(
 \omega_1\). Applying the same argument to the action of the other
 generators of \(H_1(\Sigma_g)\), we find that 
\begin{equation*}
\Omega = \prod _{i=1}^g (1 + a_i \wedge b_i \otimes u)
\end{equation*}
generates \(A = \{x \in  HFK^\infty (B) \ts | \ts \gamma \cdot x = 0 \ \
\forall \ \gamma \in H_1(\Sigma_g) \} \) as a \(\Z[u,u^{-1}] \) module. 
For future reference we note that
 the Alexander grading of \( \Omega \) has the same
parity as \(g\). 

When \(n\gg0\), the knot Floer homology tell us that 
\( \hfp (B_{-n}, \vspi_k') \cong H_*(C_k)\), where \(C_k\) is the quotient
complex of \( HFK^\infty (B) \) spanned by 
\begin{equation*}
 \{ \omega \otimes u^n \ts | \ts \omega \in \Lambda^{g+i}
 (H^1(\Sigma)), \ts n\geq \max \{k-i,0\} \}.
\end{equation*}
Moreover, this isomorphism respects the \(H_1\) action.
There are no differentials in this complex, so \( H_*(C_k) \cong
C_k\).

Let \(\pi(A)\) be the image of \(A\) in \(C_k\). We claim that for \(k
 \geq -g\), 
 the minimum Alexander grading of a nonzero element of
\(\pi(A)\) is 
\begin{equation*}
m_k = - g + 2 \Bigl\lceil \frac{g+k}{2} \Bigr\rceil.
\end{equation*}
 Indeed, it is not difficult to see that \(m_k\) is the minimum Alexander
grading of {\it any} element in \(C_k\) with the same parity as \(g\),
and that this grading is realized by any element in
\(\Lambda^{g+ m_k}(H^1(\Sigma_g))\). The expansion of \(\Omega \) certainly
contains terms of the form \( \omega \otimes u^n\), \( \omega \in 
\Lambda^{g+ m_k}(H^1(\Sigma_g))\), so the element of \(A\) with
Alexander grading \(m_k\) has a nontrivial image in \(C_k\). This proves
the claim. 

Since \(B\) is perfect, the absolute grading on \(\hfp(B_{-n},
\spi_k)\) coincides with the Alexander grading on \(C_k\) up to an
overall shift. We claim that  for \(k\leq0\), this shift is  \(E(n,k)\). 
Indeed, for \(k\leq 0\), the map 
\begin{equation*}
 F^+_{X_2,\vspi_k} \co \hfp
(\#^{2g}(S^1\times S^2), \vspi) \to \hfp (B_{-n}, \vspi_k')
\end{equation*} 
is induced by the quotient
map \( HFK^+(B) \to C_k \). Now the Alexander grading on
\(HFK^+(B)\) is equal to the absolute grading on \(\hfp
(\#^{2g}(S^1\times S^2))\), and \(F^+_{X_2,\vspi_k}\) shifts the absolute
grading by  \(E(n,k)\). This proves the claim, and thus the
proposition, when \(k\leq 0\). 
Finally, when \(k > 0\), the result follows from the conjugation symmetry
of \( \hfp\).
\end{proof}

\section{Invariants of lens spaces}
\label{Sec:Dedekind}

In this section, we establish the various properties of the
 Casson--Walker and \(d\)--invariants of lens spaces which were
used in section~\ref{Sec:Proof}. For the most part, the proofs involve
little more than elementary arithmetic. Our starting point is the
 recursive formula for the \(d\)--invariants of a lens space. In
 \cite{OS4} \Ozsvath and \Szabo introduce  natural maps from \(\Z\) to
 the set of \( \spinc \) structures on any \(L(p,q)\), which send an
 integer \(i\) to a \( \spinc\) structure \( \spi_i\). These maps have
 the property that \(\spi_i = \spi_j \) whenever \(i \equiv j \ts (p)
 \). With this labeling, they prove

\begin{prop}[Proposition 4.8 of \cite{OS4}] 
\label{Prop:dForm}Suppose \(p>q>0\) are relatively prime
 integers, and
\(0\leq i< p+q \). Then 
\begin{equation*}
d(L(p,q),\spi_i) = \frac{1}{4} - \frac{(2i+1-p-q)^2}{4pq} - d(L(q,p),\spi_i).
\end{equation*}
\end{prop}
\noindent (Note that our orientation convention for lens spaces 
is the opposite of the one in \cite{OS4}.) 

Together with the fact that
\(d(L(1,1),\spi_0) = d(S^3) = 0\), this relation clearly
determines \(d(L(p,q),\spi_i)\) for any values of \(p,q\), and
\(i\) such that \(p\) and \(q\) are relatively prime. 
For the remainder of this section,
 we adopt the shorthand notation \(d(p,q,i)\) to
stand for \(d(L(p,q),\spi_i)\). 

The Casson--Walker invariant also satisfies a recursive
formula. To be specific, \(\lambda (L(p,q)) \) is
given by a classical arithmetic function, known as a {\it Dedekind
  sum} \cite{Walker}, and it is well known that this
function satisfies a recursion relation \cite{RadGross}. For our
purposes, this relation can be stated as follows:
\begin{prop}
 Suppose \(p>q>0\) are relatively
  prime. Then
\begin{equation*}
\lambda(L(p,q)) = \frac{1}{4} - \frac{p^2+q^2+1}{12pq} - \lambda (L(q,p)).
\end{equation*}
\end{prop}
\noindent Again, it is clear that this formula, together with the condition \(
\lambda (L(1,1)) = \lambda (S^3) = 0 \) is sufficient to determine \(
\lambda (L(p,q)) \) for all relatively prime \(p\) and \(q\). As with
\(d\), we adopt the shorthand notation 
\( \lambda(p,q)\) for  \( \lambda(L(p,q)).\)

\subsection{Preliminaries} Our first order of business is to show that \(d\) and \( \lambda \)
are related:

\begin{proof}[Proof of Lemma~\ref{Lem:DCW}] Let
\begin{equation*}
\tilde{\lambda} (p,q) = \frac{1}{p} \sum_{i=0}^{p-1} \ts d(p,q,i).
\end{equation*}
We will show that \( \tilde{\lambda} \) satisfies the same recursion
relation as \( \lambda \). We write
\begin{equation*}
 \tilde{\lambda}(p,q) = 
\frac{1}{pq} \sum_{j=0}^{q-1} \ts \sum _{i=0}^{p-1} \ts d(p,q,i+j).
\end{equation*}
Applying the recursion formula and switching the order of summation, we get
\begin{align*}
\tilde{\lambda}(p,q) & = 
\frac{1}{pq} \sum_{j=0}^{q-1} \ts  \sum _{i=0}^{p-1} \ts \Bigl[\frac{1}{4} - 
\frac{(2(i+j)+1 - p - q)^2}{4pq}\Bigr] - 
\frac{1}{pq} \sum_{i=0}^{p-1} \ts \sum _{j=0}^{q-1} \ts d(q,p,i+j) \\
  & = \frac{1}{4} - 
\frac{1}{pq} \sum_{j=0}^{q-1} \ts \sum _{i=0}^{p-1} \ts
\frac{(2(i+j)+1 - p - q)^2}{4pq} - \tilde{\lambda}(q,p).
\end{align*}
Using standard identities (or simply asking Mathematica) one
finds  that 
\begin{equation*}
\sum_{j=0}^{q-1} \ts  \sum _{i=0}^{p-1} \ts  \frac{(2(i+j)+1 - p - q)^2}{4pq}
= \frac{p^2+q^2+1}{12}\cdot
\end{equation*}
This proves the claim.
\end{proof}

To  estimate the size of   \(\lambda(p,q)\), we will express it using continued fractions. To be precise, we
consider the Hirzebruch--Jung continued fraction expansion
\begin{equation*}
p/q = [a_1,a_2,\ldots ,a_n] =  a_1 - \cfrac{1}{a_2 - \cfrac{1}{\ldots - \cfrac{1}{a_n}}}
\qquad (a_i\geq 2)
\end{equation*}
common in the theory of lens spaces. The \(a_i\) may be found
recursively using the division algorithm:
\begin{equation*}
p_i/q_i = a_i - q_{i+1}/p_{i+1}
\end{equation*}
where \(p/q = p_1/q_1\),
 \(p_{i+1} = q_i\) and \(0 < q_{i+1} < p_{i+1} \). Then we
have
\begin{lem}
\label{Lem:CFExp}
\begin{equation*}
\lambda(p,q) = -\frac{1}{12}\Bigl(\frac{q}{p} + \frac{q'}{p} +
\sum_{i=1}^n (a_i - 3) \Bigr)
\end{equation*}
where \(1\leq q' <p\) and \(qq' \equiv 1 \ts (p)\). 
\end{lem}
The existance of such formulae is well known (see for example \cite{Barkan}, 
\cite{Walker}, \cite{KirbyMelvin}.) For the reader's convenience, we sketch an
elementary proof here.
\begin{proof}
We induct on the length \(n\) of the continued fraction expansion. If
\(n=1\), then \(q=1\) and \(p=a_1\), and we have
\begin{align*}
\lambda(p,1) & = \frac{1}{4} - \frac{p^2+2}{12p} \\
& = - \frac{1}{12}\Bigl(\frac{1}{p} +\frac{1}{p}+  p-3 \Bigr)
\end{align*}
which agrees with the stated form. 

In general, we note that \(q_2 \equiv - p_1 \ts (q_1)\) and that 
\begin{equation*}
\lambda (p,-q) = \lambda(L(p,-q)) = \lambda (L(p,q)^*) = -
\lambda(p,q)
\end{equation*}
so the recursion relation becomes
\begin{align*}
\lambda(p_1,q_1) & = \frac{1}{4} - \frac{p_1^2+q_1^2 + 1}{12p_1q_1} -
\lambda(q_1,p_1) \\
           & = - \frac{1}{12}\Bigl(\frac{p_1}{q_1} + \frac{q_1}{p_1}
+ \frac{1}{p_1q_1} - 3\Bigr) + \lambda(p_2,q_2).
\end{align*}
Applying the induction hypothesis to \( \lambda(p_2,q_2)\), we get
\begin{equation*}
\lambda (p_1,q_1) =  - \frac{1}{12}\Bigl(\frac{p_1}{q_1} + \frac{q_1}{p_1}
+ \frac{1}{p_1q_1} - 3 + \frac{q_2}{p_2} + \frac{q_2'}{p_2} + 
\sum_{i=2}^n (a_i - 3)\Bigr).
\end{equation*}
Now by definition
\begin{equation*}
\frac{p_1}{q_1} + \frac{q_2}{p_2} = a_1,
\end{equation*}
and it is not difficult to see that 
\begin{equation*}
\frac{1}{p_1q_1}+\frac{q_2'}{p_2} = \frac{q_1'}{p_1}.
\end{equation*}
Substituting these relations into the equation above, we obtain the
desired formula. 
\end{proof}

\subsection{Proof of Proposition~\ref{Prop:DBound}} 
Suppose that 
\begin{equation*}
\lambda(p,q) \leq \frac{1}{4}(\frac{p}{4}-1) + \lambda(p,1).
\end{equation*}
Substituting the formula of Lemma~\ref{Lem:CFExp} and simplifying, we
find that
\begin{equation*}
\frac{q}{p}+ \frac{q'}{p} + \sum_{i=1}^n (a_i-3)
\geq \frac{p}{4}+\frac{2}{p}.
\end{equation*}
Since the two fractions on the left-hand side are both \(<1\), this
implies that
\begin{equation*}
\sum_{i=1}^n (a_i-3) > \frac{p}{4} - 2.
\end{equation*}
Thus for \(p>100\), we see that
\begin{equation*}
S = \sum_{a_i>3} (a_i-3) \geq \sum_{i=1}^n (a_i-3) >\frac{2p}{9}\cdot
\end{equation*}
We investigate the conditions which this inequality puts on the
continued fraction expansion. Our first step is to estimate  the size
of \(p\) in terms of the \(a_i\).

\begin{lem}
\label{Lem:Small}
\( p_i \geq p_{i+1}(a_i-1) \).
\end{lem}

\proof
We have 
$$
p_i = p_{i+1}(a_i- \frac{q_{i+1}}{p_{i+1}})
\geq p_{i+1}(a_i-1).
\eqno{\qed}$$

\begin{cor}
\label{Lem:Product}
\({\prod_{i=1}^{n}(a_i-1) \leq  p}.\)
 Moreover, the
inequality still holds if one (but not both) of
 the factors \(a_1-1\) and \(a_n-1\) is replaced by \(a_1\) or
 \(a_n\), respectively. 
\end{cor}
\begin{proof}
An obvious induction. The case where \(a_1-1\) is replaced by \(a_1\) follows from the fact that the continued fractions \([a_1,a_2,\ldots,a_n]\) and 
\([a_n,a_{n-1},\ldots,a_1]\) have the same numerator.
\end{proof}

\begin{lem}
\label{Lem:OneBig}
If \(S > 2p/9\),
then at most two \(a_i\) are greater than \(3\). 
\end{lem}

\begin{proof}  Suppose more than one of
  the \(a_i\) is \(>3\). If there are \(m\) such terms, it is clear
  that they must all be less than \(p/3^{m-1}\). Then
\begin{equation*}
\sum_{i=1}^n (a_i-3) \leq m\cdot \frac{p}{3^{m-1}} \leq
\frac{4p}{27}
\end{equation*}
if \(m \geq 4\). Now supposing that \(m=3\), we try to maximize
\(a_1+a_2+a_3 - 9\) subject to the constraints
\((a_1-1)(a_2-1)(a_3-1) \leq p\), \(a_i\geq 4\). It is not
difficult to see that the maximum is \(\frac{p}{9}-1\) (attained when
two of the three are equal to~\(4\)), so this case is ruled out as
well.
\end{proof}

\begin{lem}
\label{Lem:A}
Let \(A = \{a_i \ts | \ts a_i >2 \} \), and let \(x\) be the largest of
the \(a_i\). If  \(S > 2p/9\)
then \(A\) is equal to one of  \(\{x\}\), \(\{x,3\}\), or \(\{x,4\}\). 
\end{lem}

\begin{proof}
Clearly \(x>3\), or \(S = 0\). 
Suppose that two of the \(a_i\), say \(x\) and \(y\), are \(>3\). If
\( y>5\), then the same sort of maximization argument used in the
previous lemma shows that \(S \leq p/5  \). If \(y=4\) or \(5\), and
 one of the other \(a_i=3\) in addition, then \((x-1)(y-1) \leq p/2\), and it
follows that \(S \leq p/6\). Finally, suppose that \(x\) is the only
value of \(a_i > 3\). Then if three or more of the other \(a_i\) equal
\(3\), we have \(S \leq p/8\). It follows that \(A\) must be one of 
\( \{x\}\), \(\{x,3\}\), \(\{x,4\}\), \(\{x,5\}\), or \( \{x,3,3\} \).

To elimate the last two possibilities, we use the sharper version of
Corollary~\ref{Lem:Product}. For example, if \(A=\{x,3,3\}\),
then one of \(a_1\) or \(a_n\) is equal to \(2\) or \(3\). If 
it is \(2\), we must have \(x \leq p/8\), whence \(S \leq p/8 \) as
well. If it is \(3\), we get that \(S\leq p/6\). A similar argument takes
care of the case \( A = \{x,5\}\). 
\end{proof}

In all remaining cases, we have \(S < x\). To analyze these
cases, suppose \(x=a_k\), and call the continued fractions
\([a_1,a_2,\ldots,a_{k-1}]\) and \( [a_{k+1},a_{k+2},\ldots,a_n]\) the
  {\it head} and {\it tail}, respectively.

\begin{lem}
If \(x > 2p/9 \), the numerator of the head and tail must both be less
than \(5\). 
\end{lem}
\begin{proof}
In the case of the tail, this follows immediately from
Lemma~\ref{Lem:Small}. To get the same result for the head, use the
fact that a continued fraction and its inverse have the same
numerator. 
\end{proof}

Thus there are only six possibilities for  the head and tail:
 \([ \ ],[2],[3],[2,2],[4],\) and \( [2,2,2]\) (corresponding to the fractions 
\(1/1,2/1,3/1,3/2,4/1\), and \(4/3\), respectively.) Since a continued
 fraction and its inverse correspond to the same lens space, we need only
 consider one element of each such pair. Thus there are 21 possible
 head--tail combinations. It is not difficult to check
 that 14  of these 21 have \(S\leq p/6\). The
 remaining 7 possiblities are listed in table \ref{table1},
which shows the continued
fraction expansion, the associated fraction \(p/q\), and the
difference \( \Delta = 12(\lambda(p,q) - \lambda(p,1)).\)

\begin{table}\begin{equation*}
\begin{array}{|l|c|c|c|}
\hline
\rule{0pt}{14pt} [a_1,a_2,\ldots,a_n] & p & q & \Delta  \\
\hline
\hline
\ \ [x] & x & 1 & 0 \\
\hline
\ \ [x,2]  & 2x-1 & 2 & x - x/p \\
\hline
\ \ [x,3] & 3x-1 & 3 & 2x-1 - (x+1)/p \\
\hline
\ \ [x,4] & 4x-1 & 4 & 3x-2 - (x+2)/p \\
\hline
\ \ [x,2,2] & 3x-2 & 3 &  2x - 2x/p \\
\hline
\ \ [x,2,2,2] & 4x-3 & 4 & 3x - 3x/p  \\
\hline
\ \ [2,x,2] & 4x-4 & 2x-1 & 3x \\
\hline
\end{array}
\end{equation*}\nocolon\caption{}\label{table1}\end{table}

The only expansions which have \(\Delta \leq 3 (\frac{p}{4}-1) \) are
those corresponding to \(L(p,1)\), \(L(p,2)\), and \(L(p,3)\). It
follows that 
the proposition holds for all values of \(p>100\). Using a computer,
it is elementary to check that it holds for all values of \(p \leq
100\) as well. This concludes the proof of the proposition. \qed

\subsection{Proof of Proposition~\ref{Prop:23}} 
To rule out the exceptional cases \(L(p,2)\) and \(L(p,3)\), 
we return to considering the \(d\)--invariants. 
 Suppose that \(L(p,2)\) is \(K_{-p}\) for some knot \(K\). Then 
\begin{equation*}
\{d(p,2,i) \ts | \ts i \in \Z/p\} = 
\{d(p,1,i) + 2n_i \ts | \ts i \in \Z/p\}
\end{equation*}
where the \(n_i\) are non-negative integers.
Since the continued fraction expansions of \(p/1\)
and \(p/2\) are short, it is easy to use the formula of
Proposition~\ref{Prop:dForm} to work out \(d(p,1,i)\) and
\(d(p,2,i)\). We get 
\begin{equation*}
d(p,1,i) = \frac{1}{4} \Bigl(1-\frac{(2i-p)^2}{p}\Bigr)
\end{equation*}
and
\begin{equation*}
d(p,2,i) = \begin{cases}
\vrule width 0pt depth 12pt \frac{1}{4}  \Bigl(2-\frac{(2i-p-1)^2}{2p}\Bigr) & \text
{if \(i\) is even,} \\
\frac{1}{4}  \Bigl(-\frac{(2i-p-1)^2}{2p}\Bigr) & \text
{if \(i\) is odd.}
\end{cases} 
\end{equation*}
We consider the largest values attained by \(L(p,1,i)\). Since we are working with
  \(L(p,2)\), \(p\) is odd, and 
the largest possible value of \(L(p,1,i) \) is
\((1-\frac{1}{p})/4\). By hypothesis, this is equal to
\(d(p,2,i)-2n_i\) for some value of \(i\). The only way this can
happen is if \(i\) is even and \(n_i=0\). In this case, we get
\begin{equation*}
1-\frac{1}{p} = 2 - \frac{(2i-p-1)^2}{2p}\cdot
\end{equation*}
After some simplification, this becomes
\begin{equation*}
(2i-p-1)^2=2p+2
\end{equation*}
so \(2p+2\) is a perfect square. We now apply the same argument to
the second largest value of \(d(p,1,i)\), which is
\((1-\frac{9}{p})/4\). Assuming \(p >9\), we again see that this must
be equal to \(d(p,2,i')\), where \(i'\) is even. Thus we have
\begin{equation*}
1-\frac{9}{p} = 2 - \frac{(2i'-p-1)^2}{2p}
\end{equation*}
which reduces to 
\begin{equation*}
(2i'-p-1)^2=2p+18.
\end{equation*}
So \(2p+2\) and \(2p+18\) are even perfect squares differing by
\(16\). But this is impossible, since the only pair of squares with
this property  is \(0\) and \(16\). Thus we need only consider those
values of  \(p\) which are \(\leq  9\). 
Consulting the list at the end of \cite{OS4}, we see
that \(L(7,2) = L(7,4) \) is the only case in which \(L(p,2)\) can be
realized as \(-p\) surgery on a knot. 

The  proof for \(L(p,3)\) is similar in spirit, but
involves a larger number of cases. We have 
\begin{equation*}
d(L(p,3,i)) = \frac{1}{4}\Bigl(1-\frac{(2i-p-2)^2}{3p}\Bigr) - d(3,p,i).
\end{equation*}
Suppose \(p \equiv 1 \ts (6) \). Then \(d(3,p,i) = -1/2 \) if \(i
\equiv 0 \ts (3)\) and \(1/6\) otherwise. We need
\begin{equation*}
\frac{1}{4} (1-\frac{1}{p}) = d(3,p,i) +2 n_i  
\qquad \text{and}\qquad \frac{1}{4}(1-\frac{9}{p}) =
d (p,3,i') + 2 n_{i'}.
\end{equation*}
If \(p > 14 \), this can only occur if \(i\) and \(i'\) are divisible
by \(3\) and \(n_i = n_i' = 0 \). In this case, we find
\begin{equation*}
(2i-p-2)^2 = 6p+3 \qquad \qquad (2i'-p-2)^2 = 6p+27
\end{equation*}
so we are looking for a pair of perfect squares which differ by
24. Again, it is easy to see that there is no such pair with the
right values\(\mod 6\). Among the possible values of \(p \equiv 1 \ts
(6) \) less than 14, the table in \cite{OS4} shows that \(L(13,3)\cong
L(13,9) \) can be a lens space surgery, while \( L(7,3)\) cannot. 

If \(p \equiv 5 \ts (6) \), a similar analysis leads to the equations
\begin{equation*}
(2i-p-2)^2 = 2p+3 \qquad \qquad (2i'-p-2)^2 = 2p+27
\end{equation*}
which actually has a solution when \(p=11\). Finally,
the remaining cases\break \( p
\equiv 2,4 \ts (6) \) do not admit any solutions. 
\qed

\end{document}